\documentclass[12pt]{amsart}
\usepackage[a4paper]{geometry}
 \newtheorem{thm}{Theorem}[section]
 
 \newtheorem{lem}[thm]{Lemma}
 
 \theoremstyle{definition}
 
 \theoremstyle{remark}
 
 \numberwithin{equation}{section}
\newcommand{\grad}{\mathop{\mathrm{grad}}}
\newcommand{\trace}{\mathop{\mathrm{trace}}}
\newcommand{\dive}{\mathop{\mathrm{div}}}

\begin{document}
\title {\textbf{Biconservative Lorentz hypersurfaces in $\mathbb{E}_{1}^{\lowercase{n}+1}$ with complex eigenvalues}}
\author{Ram Shankar Gupta, A. Sharfuddin}
 \maketitle
\begin{abstract}
Our paper is an attempt to classify biconservative submanifolds and
biharmonic submanifolds. We prove that every biconservative Lorentz
hypersurface $M_{1}^{n}$ in $\mathbb{E}_{1}^{n+1}$ having complex
eigenvalues has constant mean curvature. Moreover, every biharmonic
Lorentz hypersurface $M_{1}^{n}$ having complex
eigenvalues in $\mathbb{E}_{1}^{n+1}$ must be minimal. \\
\\
\textbf{AMS2010 MSC Codes:} 53D12, 53C40, 53C42\\
\textbf{Key Words:} Lorentz hypersurfaces, Biconservative
submanifolds, Biharmonic submanifolds, Mean curvature vector.
\end{abstract}

\maketitle
\section{\textbf{Introduction}}
The classification of constant mean curvature (CMC) hypersurfaces
play an important role in relativity theory \cite{1,2} and such type
of hypersurfaces are associated with the problem of eigenvalues of
the shape operator or differential equations arises from Laplacian
operator.

In 1964, Eells and Sampson \cite{5} introduced the notion of
poly-harmonic maps as a natural generalization of the well-known
harmonic maps.
 Thus, while harmonic maps between Riemannian manifolds $\phi: (M, g)\rightarrow(N, h)$
 are critical points of the energy functional $E(\phi)=\frac{1}{2}\int_{M}|d\phi|^{2}v_{g}$,
 the biharmonic maps are critical points of the bienergy functional
 $E_{2}(\phi)=\frac{1}{2}\int_{M}|\tau(\phi)|^{2}v_{g}$, where $\tau=\trace\nabla d\phi$ is the tension field of
 $\phi$.

In 1924, Hilbert pointed that the stress-energy tensor associated to
a functional $E$, is a conservative symmetric 2-covariant tensor $S$
at the critical points of $E$, i.e. $\dive S = 0$ \cite{3}. For the
bienergy functional $E_2$, Jiang defined the stress-bienergy tensor
$S_2$ and proved that it satisfies $\dive S_2 =
-\langle\tau_2(\phi), d\phi\rangle$ \cite{4}. Thus, if $\phi$ is
biharmonic, then $\dive S_2 = 0$. For biharmonic submanifolds, from
the above relation, we see that $\dive S_2 = 0$ if and only if the
tangent part of the bitension field vanishes. In particular, an
isometric immersion $\phi: (M, g)\rightarrow(N, h)$ is called
biconservative if $\dive S_2 = 0$.

 The biconservative submanifolds were studied and
classified in $\mathbb{E}^4$ by Hasanis and Vlachos \cite{6} in
which the biconservative hypersurfaces were called
$H$-hypersurfaces. In \cite{k2}, the complete classification of
H-hypersurfaces with three distinct curvatures in Euclidean space of
arbitrary dimension was obtained and some explicit example was
given. Upadhyay and  Turgay, classified biconservative hypersurfaces
in  $E^5_2$ with diagonal shape operator having three distinct
principal
 curvatures \cite{13}. Further, they have constructed the example of
 biconservative hypersurfaces with four distinct principal
 curvatures. Recently in \cite{k1}, it was proved that every biconservative
Lorentz hypersurface in
 $E^{n+1}_1$ with complex eigenvalues having at most five distinct principal
 curvatures has constant mean curvature. Further, it was proved that biconservative Lorentz hypersurface with
constant length of second fundamental form and whose shape operator
has complex eigenvalues with six distinct principal curvatures has
constant mean curvature \cite{k1}. For more work on biconservative
hypersurfaces in pseudo-Euclidean spaces (please see references in
\cite{k1,13}). For work on biharmonic submanifolds (please see
\cite{biha7,biha15}, and references therein).

   In
this paper, we study biconservative Lorentz hypersurfaces in
$\mathbb{E}_1^{n+1}$  whose shape operator has complex eigenvalues.
 The shape operator of Lorentz hypersurfaces with complex
eigenvalues takes the form \cite{b6,15}
\begin{equation}\label{e:1}
 \mathcal{A} = \left(
                            \begin{array}{ccccccccc}
                               \lambda & \mu & \\
                              -\mu & \lambda & & \\
                              & &      D_{n-2}&\\
                              \end{array}
                          \right), \hspace{.5 cm}
 \end{equation}
with respect to a suitable orthonormal base field of the tangent
bundle  $\{e_{1}, e_{2},..., e_{n}\}$ of $T_{p} M_{1}^{n}$, which
satisfies
\begin{equation}\label{e:2}
 g(e_{1},e_{1})=-1, \quad g(e_{i},e_{i})=1, \quad
 i=2, 3,\dots,n,
 \end{equation}
and
 \begin{equation}\label{e:3}
  g(e_{i},e_{j})=0, \quad\forall\quad i, j=1, 2,\dots,n,\quad i \neq
  j,
\end{equation}
where $D_{n-2}=diag\{\lambda_3,\lambda_4,\dots,\lambda_n\}$ and
$\mu\neq 0$.

 We prove

  \begin{thm}
Let $M_{1}^{n}$ be a biconservative Lorentz hypersurface  in
$\mathbb{E}_{1}^{n+1}$ with complex eigenvalues. Then, it has
constant mean curvature.
\end{thm}

The submanifolds satisfying
\begin{equation}\label{e:p2}
\triangle \vec {H} = 0,
\end{equation}
is called biharmonic submanifold.

The study of biharmonic submanifolds in Euclidean spaces was
initiated by Chen in  mid 1980s. In particular, he posed the
following well-known conjecture in 1991:

  \emph{The only
biharmonic submanifolds of Euclidean spaces are the minimal ones.}

 The conjecture was later studied by many researchers and so far it is found to be true for hypersurfaces in
 Euclidean spaces \cite{r20,r13,biha15,r18,r19,6}.
 Chen's
conjecture is not true always for the submanifolds of the
semi-Euclidean spaces (see \cite{r3,r5,biha6}). However, for
hypersurfaces in semi-Euclidean spaces, Chen's conjecture is also
right (see \cite{r1,r5,biha6,r23,r7,r14}).

Since, every biconservative hypersurface is a biharmonic
hypersurface, therefore, using Theorem 1.1 and normal part
$\triangle H +  H \trace (\mathcal{A}^{2}) = 0$ of the biharmonic
equation $\triangle\vec{ H}=\vec{0}$, we find
  \begin{thm}
Every biharmonic Lorentz hypersurface $M_{1}^{n}$ in
$\mathbb{E}_{1}^{n+1}$ with complex eigenvalues must be minimal.
 \end{thm}

\section{\textbf{Preliminaries}}

   Let $(M_{1}^{n}, g)$ be an $n$-dimensional Lorentz hypersurface isometrically immersed in  $(\mathbb{E}^{n+1}_{1}, \overline g)$ and $g = \overline
g_{|M_{1}^{n}}$. We denote by $\xi$ unit normal vector to
$M_{1}^{n}$ where $\overline{g}(\xi, \xi)= 1$. A vector $X$ in
$\mathbb{E}_{1}^{n+1}$ is called
  spacelike, timelike or lightlike according as \hspace{.2cm}$\overline g(X, X)> 0, \hspace{.2cm} \overline g(X, X)<0$ or \hspace{.2cm} $\overline g(X,
  X)=0$, respectively.

 The mean curvature $H$ of $M_{1}^{n}$ is given by
\begin{equation}\label{e:c4}
H = \frac{1}{n} \trace \mathcal{A},
\end{equation}
where $\mathcal{A}$ is the shape operator  of $M_{1}^{n}$.

 Let $\nabla$ denotes the Levi-Civita connection on $M_{1}^{n}$. Then, the Gauss and Codazzi equations are given by
\begin{equation}\label{e:c5}
R(X, Y)Z = g(\mathcal{A} Y, Z) \mathcal{A} X - g(\mathcal{A} X, Z)
\mathcal{A} Y,
\end{equation}
\begin{equation}\label{e:c6}
(\nabla_{X}\mathcal{A})Y = (\nabla_{Y}\mathcal{A})X,
\end{equation}
respectively, where $R$ is the curvature tensor and
\begin{equation}\label{e:c7}
(\nabla_{X}\mathcal{A})Y = \nabla_{X}\mathcal{A} Y-
\mathcal{A}(\nabla_{X}Y)
\end{equation}
for all $ X, Y, Z \in \Gamma(TM_{1}^{n})$.

 The submanifolds satisfying $\triangle \vec {H} = 0$,
is called biharmonic submanifold \cite{biha7}. The biharmonic
equation can be decomposed into its normal and tangent part. Then,
the submanifolds satisfying the tangential part of the biharmonic
equation is called biconservative. Therefore, the biconservative
Lorentz hypersurfaces $M_{1}^{n}$  in $\mathbb{E}^{n+1}_1$
characterized by
\begin{equation}\label{e:c9}
\mathcal{A} (\grad H)+ \frac{n}{2}  H \grad H = 0.
\end{equation}

\section{\textbf{Biconservative Lorentz hypersurfaces  in $\mathbb{E}^{n+1}_1$}}

 In this section, we study biconservative Lorentz hypersurface in
$\mathbb{E}_1^{n+1}$ with complex eigenvalues. Since every
hypersurface with constant mean curvature is always biconservative,
therefore, we assume that the mean curvature is not constant and
$\grad H\neq0$. Assuming non-constant mean curvature implies the
existence of an open connected subset $U$ of $M^n_1$ with grad$_{x}H
\neq 0$, for all $x\in U$. From (\ref{e:c9}), it is easy to see that
$\grad H$ is an eigenvector of the shape operator $\mathcal{A}$ with
the corresponding principal curvature $-\frac{n}{2}H$. Therefore,
there do not exist biconservative Lorentz hypersurface $M^n_1$ in
$\mathbb{E}^{n+1}_1$ with two distinct principal curvatures of
non-constant mean curvature with complex
 eigenvalues.  Without losing generality, we choose
$e_{n}$ in the direction of $\grad H$. Then, the shape operator
$\mathcal{A}$ of hypersurfaces $M_1^n$ in $\mathbb{E}^{n+1}_1$ will
take the following form with respect to a suitable orthonormal frame
$\{e_{1}, e_{2},\dots, e_{n}\}$
\begin{equation}\label{e:b1} \mathcal{A}_{H}e_1=\lambda e_1-\mu e_2,\quad \mathcal{A}_{H}e_2=\mu e_1+\lambda e_2, \quad \mathcal{A}_{H}e_i=\lambda_ie_i, \quad i\in
D,
\end{equation}
where $D=\{3, 4,\dots, n\}$.

Also, we denote  the following sets by
$$ A=\{1, 2,\dots,n\},\quad  B=\{1, 2, \dots, n-1\}, \quad C=\{3, 4,\dots, n-1\}.$$

 The $\grad H$ can be expressed as
\begin{equation}\label{e:b2}
\grad H =\sum_{i=1}^{n} e_{i}(H)e_{i}.
\end{equation}

 As we have taken $e_{n}$ parallel to $\grad H$, consequently
\begin{equation}\label{e:b3}
e_{n}(H)\neq 0, \quad e_{i}(H)= 0, \quad i\in B.
\end{equation}

 We express
\begin{equation}\label{e:b4}
\nabla_{e_{i}}e_{j}=\sum_{m=1}^{n}\omega_{ij}^{m} e_{m}, \quad i, j
\in A.
\end{equation}

Differentiating (\ref{e:2}) and (\ref{e:3}) with respect to $e_k$
and using (\ref{e:b4}),
 we obtain
\begin{equation}\label{e:b5}
\omega_{ki}^{i}=0, \quad \omega_{ki}^{j}+ \omega_{kj}^{i} =0,
\end{equation}
for $i \neq j$ and $i, j, k \in A$.

Using (\ref{e:c4}) and (\ref{e:b1}), we obtain that
 \begin{equation}\label{e:q21}
2\lambda+\sum_{j=3}^{n-1}\lambda_j= \frac{3n H}{2}=-3\lambda_n.
\end{equation}

Now, we have

\begin{lem}
 Let $M^n_1$ be a biconservative Lorentz hypersurface  in
 $\mathbb{E}^{n+1}_1$ having the shape operator given by $(\ref{e:b1})$ with respect to a suitable orthonormal frame  $\{e_{1},
e_{2},\dots, e_{n}\}$. Then, \begin{equation}\label{e:b6}
e_{i}(\lambda_{j})=
(\lambda_{i}-\lambda_{j})\omega_{ji}^{j}=(\lambda_{j}-\lambda_{i})\omega_{jj}^{i},
\end{equation}
\begin{equation}\label{e:b7}
(\lambda_{i}-\lambda_{j})\omega_{ki}^{j}=
(\lambda_{k}-\lambda_{j})\omega_{ik}^{j},
\end{equation}
\begin{equation}\label{e:k12}
(\lambda_j-\lambda)\omega_{ij}^{1}-\mu \omega_{ij}^{2}=
(\lambda_i-\lambda)\omega_{ji}^{1}-\mu \omega_{ji}^{2}
\end{equation}
\begin{equation}\label{e:k13}
(\lambda_j-\lambda)\omega_{ij}^{2}+\mu\omega_{ij}^{1}=
(\lambda_i-\lambda)\omega_{ji}^{2}+\mu\omega_{ji}^{1},\end{equation}
\begin{equation}\label{e:k15}
(\lambda_{i}-\lambda_{j})\omega_{1i}^{j} =-\mu\omega_{i2}^{j}+
(\lambda-\lambda_j)\omega_{i1}^{j},\end{equation}
\begin{equation}\label{e:q31}
(\lambda_{j}-\lambda_{i})\omega_{1j}^{i} =-\mu\omega_{j2}^{i}+
(\lambda-\lambda_i)\omega_{j1}^{i},\end{equation}
\begin{equation}\label{e:k16}
(\lambda_i-\lambda)\omega_{1i}^{1}-\mu\omega_{1i}^{2}
=e_{i}(\lambda),\end{equation}
\begin{equation}\label{e:k17}
(\lambda_i-\lambda)\omega_{1i}^{2}+\mu\omega_{1i}^{1}
=-e_i(\mu),\end{equation}
\begin{equation}\label{e:k18}
e_{1}(\lambda_{i})=-\mu\omega_{i2}^{i}+(\lambda-\lambda_i)\omega_{i1}^{i},\end{equation}
\begin{equation}\label{e:k19}
(\lambda_{i}-\lambda_{j})\omega_{2i}^{j} =\mu\omega_{i1}^{j}+
(\lambda-\lambda_j)\omega_{i2}^{j},\end{equation}
\begin{equation}\label{e:q32}
(\lambda_{j}-\lambda_{i})\omega_{2j}^{i} =\mu\omega_{j1}^{i}+
(\lambda-\lambda_i)\omega_{j2}^{i},\end{equation}
\begin{equation}\label{e:k21}
(\lambda_i-\lambda)\omega_{2i}^{1}-\mu\omega_{2i}^{2}
=e_i(\mu),\end{equation}
\begin{equation}\label{e:k20}
(\lambda_i-\lambda)\omega_{2i}^{2}+\mu\omega_{2i}^{1}
=e_{i}(\lambda),\end{equation}
\begin{equation}\label{e:k22}
e_{2}(\lambda_{i})=\mu\omega_{i1}^{i}+(\lambda-\lambda_i)\omega_{i2}^{i},\end{equation}
\begin{equation}\label{e:k25}
e_{1}(\mu)=e_{2}(\lambda),\end{equation}
\begin{equation}\label{e:k26}
e_{1}(\lambda)=-e_{2}(\mu),\end{equation}
\begin{equation}\label{e:k27}
(\lambda-\lambda_i)\omega_{12}^{i}+\mu\omega_{11}^{i}
=(\lambda-\lambda_i)\omega_{21}^{i}-\mu\omega_{22}^{i},\end{equation}
for distinct $i, j, k \in D$ such that
$\lambda_{k}\neq\lambda_{j}\neq \lambda_i$.
\end{lem}
\emph{\textbf{Proof.}} Taking $X=e_{i}, Y=e_{j}$ in (\ref{e:c7}) and
using (\ref{e:b1}) and (\ref{e:b4}), we get
$$(\nabla_{e_{i}}\mathcal{A})e_{j}=e_{i}(\lambda_{j})e_{j}+\lambda_j(\omega_{ij}^{1}e_{1}+\omega_{ij}^{2}e_{2})-(\omega_{ij}^{1}\mathcal{A}e_{1}+\omega_{ij}^{2}\mathcal{A}e_{2})
+\sum_{k=3}^{n}(\lambda_{j}-\lambda_{k})\omega_{ij}^{k}e_{k} , \quad
i, j\in D.$$

Putting the value of $(\nabla_{e_{i}}\mathcal{A})e_{j}$ in
(\ref{e:c6}), we find
\begin{equation}\label{e:k11}
\begin{array}{lcl}e_{i}(\lambda_{j})e_{j}+\lambda_j(\omega_{ij}^{1}e_{1}+\omega_{ij}^{2}e_{2})-(\omega_{ij}^{1}\mathcal{A}e_{1}+\omega_{ij}^{2}\mathcal{A}e_{2})
+\sum_{k=3}^{n}(\lambda_{j}-\lambda_{k})\omega_{ij}^{k}e_{k}
\\=e_{j}(\lambda_{i})e_{i}+\lambda_i(\omega_{ji}^{1}e_{1}+\omega_{ji}^{2}e_{2})-(\omega_{ji}^{1}\mathcal{A}e_{1}
+\omega_{ji}^{2}\mathcal{A}e_{2})+\sum_{k=3}^{n}(\lambda_{i}-\lambda_{k})\omega_{ji}^{k}e_{k}
,\end{array}\end{equation} whereby for $i\neq j=k$ and $i\neq j \neq
k$, we obtain (\ref{e:b6}) and (\ref{e:b7}), respectively. Moreover,
using (\ref{e:b1}) in (\ref{e:k11}) and comparing the coefficients
of $e_1$ and $e_2$, we find (\ref{e:k12}) and (\ref{e:k13}),
respectively.

Next, using (\ref{e:b1}) and (\ref{e:c7}) in
$(\nabla_{e_{1}}\mathcal{A})e_{i}=(\nabla_{e_{i}}\mathcal{A})e_{1}$,
for $i\in D$, we obtain
\begin{equation}\label{e:k14}
\begin{array}{lcl}e_{1}(\lambda_{i})e_{i}+\lambda_i(\omega_{1i}^{1}e_{1}+\omega_{1i}^{2}e_{2})-(\omega_{1i}^{1}\mathcal{A}e_{1}+\omega_{1i}^{2}\mathcal{A}e_{2})
+\sum_{j=3}^{n}(\lambda_{i}-\lambda_{j})\omega_{1i}^{j}e_{j}
=\\(e_{i}(\lambda)-\mu\omega_{i2}^{1})e_1+(\lambda\omega_{i1}^{2}-e_i(\mu))e_2
-\omega_{i1}^{2}\mathcal{A}e_{2}+\sum_{j=3}^{n}\Big((\lambda-\lambda_j)\omega_{i1}^{j}-\mu\omega_{i2}^{j}\Big)e_{j}
,\end{array}\end{equation} whereby for $i\neq j$, we get
(\ref{e:k15}). Further, comparing the coefficients of $e_1$, $e_2$,
and $e_i$ and using (\ref{e:b5}), we have (\ref{e:k16}),
(\ref{e:k17}) and (\ref{e:k18}), respectively.

Also, using (\ref{e:b1}) and (\ref{e:c7}) in
$((\nabla_{e_{1}}\mathcal{A})e_{j},e_i)=g((\nabla_{e_{j}}\mathcal{A})e_{1},e_i)$,
gives (\ref{e:q31}).

Similarly, using (\ref{e:b1}) and (\ref{e:c7}) in
$(\nabla_{e_{2}}\mathcal{A})e_{i}=(\nabla_{e_{i}}\mathcal{A})e_{2}$,
for $i\in D$, we get
\begin{equation}\label{e:k23}
\begin{array}{lcl}e_{2}(\lambda_{i})e_{i}+\lambda_i(\omega_{2i}^{1}e_{1}+\omega_{2i}^{2}e_{2})-(\omega_{2i}^{1}\mathcal{A}e_{1}+\omega_{2i}^{2}\mathcal{A}e_{2})
+\sum_{j=3}^{n}(\lambda_{i}-\lambda_{j})\omega_{2i}^{j}e_{j}
=\\(e_{i}(\lambda)+\mu\omega_{i1}^{2})e_2+(\lambda\omega_{i2}^{1}+e_i(\mu))e_1
-\omega_{i2}^{1}\mathcal{A}e_{1}+\sum_{j=3}^{n}\Big((\lambda-\lambda_j)\omega_{i2}^{j}+\mu\omega_{i1}^{j}\Big)e_{j}
,\end{array}\end{equation} whereby for $i\neq j$, we get
(\ref{e:k19}). Further, comparing the coefficients of $e_1$, $e_2$,
and $e_i$ and using (\ref{e:b5}), we have (\ref{e:k20}),
(\ref{e:k21}) and (\ref{e:k22}), respectively.

Also, using (\ref{e:b1}) and (\ref{e:c7}) in
$((\nabla_{e_{2}}\mathcal{A})e_{j},e_i)=g((\nabla_{e_{j}}\mathcal{A})e_{2},e_i)$,
gives (\ref{e:q32}).

Now, using (\ref{e:b1}) and (\ref{e:c7}) in
$(\nabla_{e_{1}}\mathcal{A})e_{2}=(\nabla_{e_{2}}\mathcal{A})e_{1}$,
for $i\in D$, we obtain
\begin{equation}\label{e:k24}
\begin{array}{lcl}e_{1}(\mu)e_{1}+e_{1}(\lambda)e_{2}
+\sum_{i=3}^n\Big((\lambda-\lambda_i)\omega_{12}^{i}+\mu\omega_{11}^{i}\Big)e_{i}
\\=e_{2}(\lambda)e_{1}-e_{2}(\mu)e_{2}+\sum_{i=3}^n\Big((\lambda-\lambda_i)\omega_{21}^{i}-\mu\omega_{22}^{i}\Big)e_{i}
,\end{array}\end{equation} whereby  comparing the coefficients of
$e_1$, $e_2$, and $e_i$ and using (\ref{e:b5}), we have
(\ref{e:k25}), (\ref{e:k26}) and (\ref{e:k27}), respectively. This
completes the proof of the Lemma.\\

Next, we have
\begin{lem}
 Let $M^n_1$ be a biconservative Lorentz hypersurface  in
 $\mathbb{E}^{n+1}_1$ having the shape operator given by $(\ref{e:b1})$ with respect to a suitable orthonormal frame  $\{e_{1},
e_{2},\dots, e_{n}\}$. Then, \begin{equation}\label{e:k4}
\lambda_n\neq\lambda_k, \quad\forall\quad k\in C.
\end{equation}
\end{lem}
\emph{\textbf{Proof.}} Let $\lambda_n=\lambda_k$ for $k\in C$, then
taking $i=n$ and $j=k$ in (\ref{e:b6}), we get $$e_n(\lambda_k)=0
\quad\mbox{or}\quad e_n(H)=0, \quad\mbox{as}\quad
\lambda_n=-\frac{nH}{2},$$ which contradicts (\ref{e:b3}). Whereby
completing the proof of Lemma.\\

Using (\ref{e:b3}), (\ref{e:b4}) and the fact that $[e_{i}
\hspace{.1 cm}
e_{j}](H)=0=\nabla_{e_{i}}e_{j}(H)-\nabla_{e_{j}}e_{i}(H)=\omega_{ij}^{n}e_{n}(H)-\omega_{ji}^{n}e_{n}(H),$
for $i\neq j$, we find
\begin{equation}\label{e:b9}
\omega_{ij}^{n}=\omega_{ji}^{n},  \quad i, j \in B.
\end{equation}

\begin{lem}
 Let $M^n_1$ be a biconservative Lorentz hypersurface  in
 $\mathbb{E}^{n+1}_1$  having the shape operator given by $(\ref{e:b1})$ with respect to a suitable orthonormal frame  $\{e_{1},
e_{2},\dots, e_{n}\}$. Then, \begin{equation}\label{e:k29}
\omega_{nn}^{i}=0, \quad\forall\quad i\in A.
\end{equation}
\end{lem}
\emph{\textbf{Proof.}}  Putting $i\neq n, j = n$ in (\ref{e:b6}) and
using (\ref{e:b3}) and (\ref{e:b5}), we
 find
\begin{equation}\label{e:b13}
\omega_{nn}^{i}= 0, \quad  i\in D.
\end{equation}

Taking $i=n$ in (\ref{e:k18}) and (\ref{e:k22}) and using
(\ref{e:b3}) and (\ref{e:b5}), we find
\begin{equation}\begin{array}{lcl}\label{e:k30}
\omega_{nn}^{1}=\omega_{nn}^{2}=0.
\end{array}\end{equation}

Combining (\ref{e:b13}) and (\ref{e:k30}), we get (\ref{e:k29}).

\begin{lem}
 Let $M^n_1$ be a biconservative Lorentz hypersurface  in
 $\mathbb{E}^{n+1}_1$  having the shape operator given by $(\ref{e:b1})$ with respect to a suitable orthonormal frame  $\{e_{1},
e_{2},\dots, e_{n}\}$. Then, \begin{equation}\label{e:k31}
\omega_{22}^{n}=\omega_{11}^{n}=\omega_{12}^{n}=\omega_{21}^{n}=0.
\end{equation}
\end{lem}
\emph{\textbf{Proof.}}  Taking $i=n$ in (\ref{e:k27}) and using
(\ref{e:b9}), we get
\begin{equation}\begin{array}{lcl}\label{e:k28}
\omega_{11}^{n}=-\omega_{22}^{n}.
\end{array}\end{equation}

Taking $i=n$ in (\ref{e:k16}), (\ref{e:k17}), (\ref{e:k21}),
(\ref{e:k20})  and using (\ref{e:b5}), (\ref{e:b9}), and
(\ref{e:k28}), we find
\begin{equation}\begin{array}{lcl}\label{e:k32}
-(\lambda_n-\lambda)\omega_{12}^n+\mu \omega_{22}^n=0,\quad
(\lambda_n-\lambda)\omega_{22}^n+\mu \omega_{12}^n=0.
\end{array}\end{equation}

Solving (\ref{e:k32}), we get
\begin{equation}\begin{array}{lcl}\label{e:k33}
\omega_{12}^n=\omega_{22}^n=0.
\end{array}\end{equation}

Using (\ref{e:k33}), (\ref{e:k28}) and (\ref{e:b9}), we get
(\ref{e:k31}).

\begin{lem}
 Let $M^n_1$ be a biconservative Lorentz hypersurface  in
 $\mathbb{E}^{n+1}_1$  having the shape operator given by $(\ref{e:b1})$ with respect to a suitable orthonormal frame  $\{e_{1},
e_{2},\dots, e_{n}\}$. Then, \begin{equation}\label{e:k34}
\omega_{ij}^{1}=\omega_{ij}^{2}=\omega_{1i}^{j}=\omega_{2i}^{j}=0,\quad
i\neq j, \quad i, j\in D.
\end{equation}
\end{lem}
\emph{\textbf{Proof.}}  Using (\ref{e:k15}), (\ref{e:q31}) and
(\ref{e:b5}), we get
\begin{equation}\begin{array}{lcl}\label{e:k35}
\mu\omega_{ij}^{2}-(\lambda-\lambda_j)\omega_{ij}^1=\mu\omega_{ji}^{2}-(\lambda-\lambda_i)\omega_{ji}^1.
\end{array}\end{equation}

Similarly, using (\ref{e:k19}),
 (\ref{e:q32}) and (\ref{e:b5}), we find
\begin{equation}\begin{array}{lcl}\label{e:k36}
\mu\omega_{ij}^{1}+(\lambda-\lambda_j)\omega_{ij}^2=\mu\omega_{ji}^{1}+(\lambda-\lambda_i)\omega_{ji}^2.
\end{array}\end{equation}

Combining (\ref{e:k13}) and  (\ref{e:k36}), we obtain
\begin{equation}\begin{array}{lcl}\label{e:k37}
\omega_{ij}^1= \omega_{ji}^1.
\end{array}\end{equation}

Combining (\ref{e:k12}) and  (\ref{e:k35}), we find
\begin{equation}\begin{array}{lcl}\label{e:k38}
(\lambda_j-\lambda)\omega_{ij}^1=(\lambda_i-\lambda) \omega_{ji}^1.
\end{array}\end{equation}

Using (\ref{e:k35}), (\ref{e:k36}), (\ref{e:k37}) and (\ref{e:k38}),
we get (\ref{e:k34}). Whereby completing the proof of the Lemma.\\

Now, we find following Lemma for covariant derivative.
\begin{lem}
 Let $M^n_1$ be a biconservative Lorentz hypersurface  in
 $\mathbb{E}^{n+1}_1$  having the shape operator given by $(\ref{e:b1})$ with respect to a suitable orthonormal frame  $\{e_{1},
e_{2},\dots, e_{n}\}$. Then,
$$ \nabla_{e_{1}}e_{1}=\sum_{m\neq 1, n}\omega_{11}^m e_{m}, \nabla_{e_{1}}e_{2}=\sum_{m\neq 2, n}\omega_{12}^m
e_{m},  \nabla_{e_{1}}e_{n}=0, \nabla_{e_{n}}e_{1}=
\omega_{n1}^{2}e_{2},$$
 $$ \nabla_{e_{2}}e_{1}=\sum_{m\neq 1, n}\omega_{21}^m e_{m}, \nabla_{e_{2}}e_{2}=\sum_{m\neq 2, n}\omega_{22}^m
 e_{m}, \nabla_{e_{2}}e_{n}=0, \nabla_{e_{n}}e_{2}= \omega_{n2}^{1}e_{1},$$
 $$\nabla_{e_{1}}e_{i}=\omega_{1i}^1
e_{1}+\omega_{1i}^2e_2, \nabla_{e_{2}}e_{i}=\omega_{2i}^1
e_{1}+\omega_{2i}^2e_2, \nabla_{e_{i}}e_{i}=\sum_{m\neq
i}\omega_{ii}^m e_{m},$$
$$ \nabla_{e_{i}}e_{1}=\omega_{i1}^2 e_{2}+\omega_{i1}^i e_{i}, \nabla_{e_{i}}e_{2}=\omega_{i2}^1
e_{1}+\omega_{i2}^i e_{i},  \nabla_{e_{i}}e_{n}=-\omega_{ii}^ne_i,
\nabla_{e_{n}}e_{n}=0, $$ for $i \in C$. Moreover, \\$(a)$ if
$M_1^n$ has all distinct principal curvatures, then
$$\nabla_{e_{n}}e_{i}=0,
\nabla_{e_{i}}e_{j}=\omega_{ij}^{i} e_{i} \quad \forall \quad i, j
\in C, \quad i\neq j,$$  $(b)$ if $M_1^n$ has $'q'$ distinct
principal curvatures $\lambda\pm\sqrt{-1}\mu,
\lambda_3,\dots,\lambda_{q-1},\lambda_n$, with multiplicities
$p_3,\dots, p_{q-1}$ of $\lambda_3,\dots,\lambda_{q-1}$,
respectively, such that $p_3+p_4+\dots+p_{q-2}+p_{q-1}=n-3$. Then
$$\nabla_{e_{n}}e_{i}=\sum_{C_{i_1}}\omega_{ni}^{m}e_{m} \quad
\forall\quad i \in C_{i_1}, \quad m\neq i,$$ $$
\nabla_{e_{i}}e_{j}=\sum_{C_{i_1}}\omega_{ij}^{m} e_{m} \quad
\forall \quad i, j \in C_{i_1}, \quad i\neq j, \quad m\neq j,$$
 where $i_1=3, 2,\dots,q-1$, and $C_3=\{3,\dots,p_3+2\}, C_4=\{p_3+3,\dots,p_3+p_4+2\},
 \dots, C_{q-1}=\{p_3+p_4+\dots+p_{q-2}+3,\dots,n-1\}$, and
$\omega_{ij}^{i}$ satisfy $(\ref{e:b5})$ and $(\ref{e:b6})$.
\end{lem}
\emph{\textbf{Proof.}} (a) Let $M_1^n$ has all distinct principal
curvatures.
 Putting  $j = n$ and $k = j$ in (\ref{e:b7}) and using (\ref{e:b9}), we
 get
\begin{equation}\begin{array}{lcl}\label{e:q41}
\omega_{ji}^{n}=\omega_{ij}^{n}=0, \quad i, j \in C, \quad i\neq j.
\end{array}\end{equation}

Putting $i = n$ and $k=i$ in (\ref{e:b7}) and
 using (\ref{e:q41}) and  (\ref{e:b5}), we  find
\begin{equation}\label{e:q42}
\omega_{ni}^{j}=\omega_{in}^{j} = 0, \quad i, j \in C, \quad i\neq
j.
\end{equation}

(b) Let $M_1^n$ has $'q'$ distinct principal curvatures. Putting $i
= n$ and $k=i$ in (\ref{e:b7}), we
 obtain
\begin{equation}\begin{array}{lcl}\label{e:b17}
\omega_{in}^{j}= 0, \quad j\neq i \quad  \mbox{and}\quad j, i \in
C_{i_1}, \quad i_1=3, \dots, q-1.
\end{array}\end{equation}

 Putting  $j = n$ and $k=j$ in (\ref{e:b7}) and using (\ref{e:b9}), we
 get
\begin{equation}\begin{array}{lcl}\label{e:b19}
\omega_{ji}^{n}=\omega_{ij}^{n}=0, \quad i \in C_{i_1},\quad j \in
C_{i_2}, \quad i_1\neq i_2, \quad i_1, i_2= 3, \dots, q-1.
\end{array}\end{equation}

Taking $i \in C_{i_1}$  in (\ref{e:b7}), we
 have
\begin{equation}\begin{array}{lcl}\label{e:b18}
\omega_{ki}^{j}= 0, \quad j\neq k\quad\mbox{and}\quad j, k \in
C_{i_2},\quad i_1\neq i_2, \quad i_1, i_2= 3, \dots, q-1.
\end{array}\end{equation}

Putting $i = n$ and $k=i$ in (\ref{e:b7}) and
 using (\ref{e:b19}) and  (\ref{e:b5}), we  find
\begin{equation}\label{e:b20}
\omega_{ni}^{j}=\omega_{in}^{j} = 0, \quad j \in C_{i_1}, \quad i\in
C_{i_2}, \quad i_1\neq i_2, \quad i_1, i_2= 3, \dots, q-1.
\end{equation}

Now, using Lemma 3.3, Lemma 3.4, Lemma 3.5 and (\ref{e:q41}),
(\ref{e:q42}), (\ref{e:b17}), (\ref{e:b19}), (\ref{e:b18}) and
(\ref{e:b20}) in (\ref{e:b4}), completes the
proof of the Lemma.\\

Next, we have

\begin{lem}
 Let $M^n_1$ be a biconservative Lorentz hypersurface  in
 $\mathbb{E}^{n+1}_1$  having the shape operator given by $(\ref{e:b1})$ with respect to a suitable orthonormal frame  $\{e_{1},
e_{2},\dots, e_{n}\}$. Then, \begin{equation}\label{e:k39}
\lambda=0.
\end{equation}
\end{lem}
\emph{\textbf{Proof.}} Evaluating $g(R(e_{n},e_{1})e_{n},e_{1})$,
using (\ref{e:c5}), (\ref{e:b1}) and Lemma 3.6, we have
\begin{equation}\label{e:k40}\begin{array}{lcl}
g(\nabla_{e_{n}}\nabla_{e_{1}}e_{n}-\nabla_{e_{1}}\nabla_{e_{n}}e_{n}-\nabla_{[e_n\hspace{.1cm}
e_{1}]}e_{n}, e_{1})\\= g(Ae_{1}, e_{n}) g(Ae_{n}, e_{1}) -
g(Ae_{n}, e_{n}) g(Ae_{1}, e_{1}),\end{array}
\end{equation}
which gives
\begin{equation}\label{e:k41}
\lambda \lambda_{n}=0.
\end{equation}

Since $\lambda_n\neq0$, therefore, from (\ref{e:k41}), we find
(\ref{e:k39}). Thus completing the proof of the Lemma.

Now, using Lemma 3.7, we find following Theorem.
\begin{thm}
 There do not exist biconservative
Lorentz hypersurface $M^n_1$ in  $\mathbb{E}^{n+1}_1$ with three
distinct principal curvatures of non-constant mean curvature with
complex
 eigenvalues.
\end{thm}
\emph{\textbf{Proof.}} Let $M_1^n$ has three distinct principal
curvatures. Then, from (\ref{e:q21}) and (\ref{e:k39}), we get
$H=0$, a contradiction. Which completes the proof of the theorem.

Next, we have
\begin{lem}
 Let $M^n_1$ be a biconservative Lorentz hypersurface  in
 $\mathbb{E}^{n+1}_1$  having the shape operator given by $(\ref{e:b1})$ with respect to a suitable orthonormal frame  $\{e_{1},
e_{2},\dots, e_{n}\}$. Then,
\begin{equation}\label{e:k47}
\omega_{11}^i=\omega_{22}^i=\omega_{12}^i=\omega_{21}^i=0,
\quad\forall\quad i\in C,
\end{equation}
and
\begin{equation}\label{e:q23}
\mu=\mbox{constant}.
\end{equation}
\end{lem}
\emph{\textbf{Proof.}} Using (\ref{e:k39}) and (\ref{e:b5}) in
(\ref{e:k16}) and (\ref{e:k20}), we find
\begin{equation}\label{e:k42}
\lambda_i \omega_{11}^i=\mu \omega_{12}^i \quad\mbox{and}\quad
\lambda_i \omega_{22}^i=-\mu \omega_{21}^i,
\end{equation}
respectively.

On the other hand, adding (\ref{e:k17}) and (\ref{e:k21}), and
therein using (\ref{e:k39}), (\ref{e:k42}) and (\ref{e:b5}), we
obtain
\begin{equation}\label{e:k43}
\omega_{11}^i= \omega_{22}^i,
\end{equation}
which together with (\ref{e:k42}) gives
\begin{equation}\label{e:k44}
\omega_{12}^i=- \omega_{21}^i.
\end{equation}

Using (\ref{e:k39}), (\ref{e:k43}) and (\ref{e:k44}) in
(\ref{e:k27}), we get
\begin{equation}\label{e:k45}
\mu\omega_{11}^i=\lambda_i \omega_{12}^i.
\end{equation}

Therefore, from (\ref{e:k42}) and (\ref{e:k45}), we obtain
\begin{equation}\label{e:k46}
(\mu^2-\lambda_i^2)\omega_{11}^i=0.
\end{equation}

We claim that $\omega_{11}^i=0$. In fact, if $\omega_{11}^i\neq0$,
then $\mu^2-\lambda_i^2=0$. Which gives $\lambda_i=\pm\mu$ for all
$i\in C$. In view of Theorem 3.8, we consider the following cases:\\

\textbf{Case I.} Let $M_1^n$ has four distinct principal curvatures.
Then, using  (\ref{e:k39}) and Lemma 3.2 in (\ref{e:q21}), we obtain
$(n-3)\lambda_i=\frac{3nH}{2}$ or $\pm (n-3)\mu=\frac{3nH}{2}$.
Which on differentiating with respect to $e_n$ gives $\pm
(n-3)e_n(\mu)=\frac{3ne_n(H)}{2}$. Also, using (\ref{e:k31}) in
(\ref{e:k17}), we find  $e_n(\mu)=0$. Therefore, we obtain
$e_n(H)=0$, a contradiction.

\textbf{Case II.} Let $M_1^n$ has five distinct principal curvatures
$\lambda\pm \sqrt{-1}\mu, \lambda_3=\mu, \lambda_4=-\mu, \lambda_n$.
Then, using  (\ref{e:k39}) and Lemma 3.2 in (\ref{e:q21}), we get
$(p_3-p_4)\mu=\frac{3nH}{2}$, where $p_3$ and $p_4$ are the
multiplicities of $\lambda_3$ and $\lambda_4$, respectively. Now,
proceeding as in Case I, we get a contradiction.

\textbf{Case III.} Let $M_1^n$ has more than five distinct principal
curvatures. Then, $\lambda_i=\pm\mu$ for all $i\in C$ gives a
contradiction to more than five distinct principal curvatures.\\

 Hence
$\omega_{11}^i=0$. Using this in (\ref{e:k43}), (\ref{e:k44}) and
(\ref{e:k45}), we find (\ref{e:k47}).

Using (\ref{e:k31}), (\ref{e:k39}) and (\ref{e:k47}) in
(\ref{e:k17}), (\ref{e:k25}) and (\ref{e:k26}), we get
\begin{equation}\label{e:k48}
e_1(\mu)=e_2(\mu)=e_i(\mu)=0, \quad\forall\quad i\in D.
\end{equation}

Hence $\mu$ is constant in all direction. This completes the proof
of the Lemma.

\begin{lem}
 Let $M^n_1$ be a biconservative Lorentz hypersurface  in
 $\mathbb{E}^{n+1}_1$  having the shape operator given by $(\ref{e:b1})$ with respect to a suitable orthonormal frame  $\{e_{1},
e_{2},\dots, e_{n}\}$.  Then, $g(R(e_{n},e_{i})e_{n},e_{i})$,
$g(R(e_{n},e_{i})e_{i},e_{1})$, $g(R(e_{n},e_{i})e_{i},e_{2})$,
$g(R(e_{i},e_{1})e_{i},e_{n})$, $g(R(e_{i},e_{2})e_{i},e_{n})$,
$g(R(e_{i},e_{2})e_{i},e_{1})$ and $g(R(e_{i},e_{1})e_{i},e_{2})$
give the following:
\begin{equation}\label{e:q15}
  e_{n}(\omega_{ii}^{n})- (\omega_{ii}^{n})^{2}=  \lambda_{n} \lambda_{i},
\end{equation}
\begin{equation}\label{e:k49}
  e_{n}(\omega_{ii}^{1})= \omega_{ii}^{n}\omega_{ii}^1-\omega_{n2}^{1}\omega_{ii}^2,
\end{equation}
\begin{equation}\label{e:k50}
  e_{n}(\omega_{ii}^{2})=\omega_{ii}^{n}\omega_{ii}^2-\omega_{n1}^{2}\omega_{ii}^1,
\end{equation}
\begin{equation}\label{e:k51}
  e_{1}(\omega_{ii}^{n})=\omega_{ii}^{n}\omega_{ii}^1,
\end{equation}
\begin{equation}\label{e:k52}
  e_{2}(\omega_{ii}^{n})=\omega_{ii}^{n}\omega_{ii}^2,
  \end{equation}
\begin{equation}\label{e:k53}
  e_{2}(\omega_{ii}^{1})+ \omega_{ii}^{2}(\omega_{22}^1-\omega_{ii}^1)= \mu \lambda_{i},
\end{equation}
and
\begin{equation}\label{e:k54}
  e_{1}(\omega_{ii}^{2})+ \omega_{ii}^{1}(\omega_{11}^2-\omega_{ii}^2)= \mu \lambda_{i},
\end{equation}
respectively, for all $i \in
 C$.
\end{lem}
\emph{\textbf{Proof:}} Here, we give the proof of the first two
relations (\ref{e:q15}) and (\ref{e:k49}). The proof of the other
relations can be obtained in a similar way.

Using (\ref{e:c5}) and (\ref{e:b1}), we have
\begin{equation}\label{e:p3} g(R(e_{n},e_{i})e_{n},e_{i})
= g(Ae_{i}, e_{n}) g(Ae_{n}, e_{i}) - g(Ae_{n}, e_{n}) g(Ae_{i},
e_{i})=- \lambda_n \lambda_{i},
\end{equation}
\begin{equation}\label{e:k55} g(R(e_{n},e_{i})e_{i},e_{1})
= g(Ae_{i}, e_{i}) g(Ae_{n}, e_{1}) - g(Ae_{n}, e_{i}) g(Ae_{i},
e_{1})=0,
\end{equation} \quad for all $i\in C$.

(i) Let $M_1^n$ has all the distinct principal curvatures. Then,
using Lemma 3.6, we get
\begin{equation}\label{e:q51}\begin{array}{lcl} g(R(e_{n},e_{i})e_{n},e_{i})=
g(\nabla_{e_{n}}\nabla_{e_{i}}e_{n}-\nabla_{e_{i}}\nabla_{e_{n}}e_{n}-\nabla_{[e_{n}\hspace{.1cm}
e_{i}]}e_{n},
e_{i})\\=g(\nabla_{e_{n}}(-\omega_{ii}^{n}e_{i})-\omega_{ii}^{n}\nabla_{e_{i}}
e_{n}, e_{i})= g(-e_{n}(\omega_{ii}^{n})
e_{i}+(\omega_{ii}^{n})^{2}e_{i}, e_{i})\\=
-e_{n}(\omega_{ii}^{n})+(\omega_{ii}^{n})^{2},\end{array}\end{equation}
for all $i\in C$.

Therefore, from (\ref{e:p3}) and (\ref{e:q51}), we get
(\ref{e:q15}).

Next, we know that
\begin{equation}\label{e:k56}\begin{array}{rcl} g(R(e_{n},e_{i})e_{i},e_{1})=
g(\nabla_{e_{n}}\nabla_{e_{i}}e_{i}-\nabla_{e_{i}}\nabla_{e_{n}}e_{i}-\nabla_{[e_{n}\hspace{.1cm}
e_{i}]}e_{i}, e_{1}),\end{array}\end{equation} for all $i\in C$.

Now, using Lemma 3.6 and Lemma 3.9, we have
$$\nabla_{e_{n}}\nabla_{e_{i}}e_{i}=\nabla_{e_{n}}(\sum_{m\neq
i, m=1}^{n}\omega_{ii}^{m}e_{m})=\sum_{m\neq
i,m=1}^n\big(e_{n}(\omega_{ii}^{m})e_{m}+\omega_{ii}^{m}\nabla_{e_{n}}e_m\big),$$
$$\nabla_{e_{i}}\nabla_{e_{n}}e_{i}=0, \nabla_{\nabla e_{n}e_i}e_{i}=0,$$
$$\nabla_{\nabla e_{i}e_n}e_{i}=-\omega_{ii}^{n}\nabla_{e_{i}}e_{i}=-\omega_{ii}^{n}(\sum_{m\neq i,m=1}^{n}\omega_{ii}^{m}e_m).$$

Hence, using above in (\ref{e:k56}), we get
\begin{equation}\label{e:q61}\begin{array}{rcl}
g(R(e_{n},e_{i})e_{i},e_{1})=-e_n(\omega_{ii}^{1})-\omega_{ii}^{2}\omega_{n2}^{1}+\omega_{ii}^{n}\omega_{ii}^{1}.\end{array}\end{equation}
for all $i\in C$.

Therefore, from (\ref{e:k55}) and (\ref{e:q61}), we get
(\ref{e:k49}) for $i\in C$.\\

(ii) Let $M_1^n$ has 'q' distinct principal curvatures. Then, using
Lemma 3.6, we find
\begin{center} $\begin{array}{lcl} g(R(e_{n},e_{i})e_{n},e_{i})=
g(\nabla_{e_{n}}\nabla_{e_{i}}e_{n}-\nabla_{e_{i}}\nabla_{e_{n}}e_{n}-\nabla_{[e_{n}\hspace{.1cm}
e_{i}]}e_{n},
e_{i})\\=g(\nabla_{e_{n}}(-\omega_{ii}^{n}e_{i})-\sum_{m\neq
i}^{C_3}\omega_{ni}^{m}(\nabla_{e_{m}} e_{n})
-\omega_{ii}^{n}\nabla_{e_{i}} e_{n}, e_{i})\\=
g(-e_{n}(\omega_{ii}^{n}) e_{i}-\sum_{m\neq
i}^{C_3}\omega_{ni}^{m}(\sum_{l\neq
n}^{C_3}\omega_{mn}^{l}e_l)+(\omega_{ii}^{n})^{2}e_{i},
e_{i}),\end{array}$\end{center} wherein using (\ref{e:b17}), gives
\begin{equation}\label{e:p4}\begin{array}{lcl} g(R(e_{n},e_{i})e_{n},e_{i})=
-e_{n}(\omega_{ii}^{n})+(\omega_{ii}^{n})^{2},\end{array}\end{equation}
for all $i\in C_3$.

Therefore, from (\ref{e:p3}) and (\ref{e:p4}), we get (\ref{e:q15})
for $i\in C_3$. Similarly, for all $i\in C$, we find (\ref{e:q15}).

Next, using Lemma 3.6 and Lemma 3.9 for $i\in C_3$, we have
$$\nabla_{e_{n}}\nabla_{e_{i}}e_{i}=\nabla_{e_{n}}(\sum_{m\neq
i, m=1}^{n}\omega_{ii}^{m}e_{m})=\sum_{m\neq
i,m=1}^n\big(e_{n}(\omega_{ii}^{m})e_{m}+\omega_{ii}^{m}\nabla_{e_{n}}e_m\big),$$
$$\nabla_{e_{i}}\nabla_{e_{n}}e_{i}=\nabla_{e_{i}}(\sum_{m\neq
i,m=3}^{p_3+2}\omega_{ni}^{m}e_{m})=\sum_{m\neq
i,m=3}^{p_3+2}\big(e_{i}(\omega_{ni}^{m})e_{m}+\omega_{ni}^{m}\sum_{l\neq
m,m=3}^{p_3+2}\omega_{im}^{l}e_l\big),$$
$$\nabla_{\nabla e_{n}e_i}e_{i}=\sum_{m\neq
i,m=3}^{p_3+2}\omega_{ni}^{m}\nabla_{e_{m}}e_{i}=\sum_{m\neq
i,m=3}^{p_3+2}\omega_{ni}^{m}(\sum_{l\neq
i,l=3}^{p_3+2}\omega_{mi}^{l}e_l),$$
$$\nabla_{\nabla e_{i}e_n}e_{i}=-\omega_{ii}^{n}\nabla_{e_{i}}e_{i}=-\omega_{ii}^{n}(\sum_{m\neq i,m=1}^{n}\omega_{ii}^{m}e_m).$$

Hence, using above in (\ref{e:k56}), we get
\begin{equation}\label{e:k57}\begin{array}{rcl}
g(R(e_{n},e_{i})e_{i},e_{1})=-e_n(\omega_{ii}^{1})-\omega_{ii}^{2}\omega_{n2}^{1}+\omega_{ii}^{n}\omega_{ii}^{1}.\end{array}\end{equation}
for all $i\in C_3$.

Therefore, from (\ref{e:k55}) and (\ref{e:k57}), we get
(\ref{e:k49}) for $i\in C_3$. Similarly, for all $i\in C$, we find
(\ref{e:k49}).

\section{\textbf{Proof of the theorem}}

 Using Lemma 3.6, we get
\begin{equation}\label{e:m1}
e_1e_n-e_ne_1=\nabla_{e_1}e_n-\nabla_{e_n}e_1=-\omega_{n1}^{2}e_2.\end{equation}

Operating $\omega_{ii}^n$ on both sides in (\ref{e:m1}), we find
\begin{equation}\label{e:m2}
e_1e_n(\omega_{ii}^n)-e_ne_1(\omega_{ii}^n)=-\omega_{n1}^{2}e_2(\omega_{ii}^n).\end{equation}

Using (\ref{e:q15}), (\ref{e:k49}),  (\ref{e:k51}), (\ref{e:k52}),
(\ref{e:k18}), (\ref{e:b3}) and Lemma 3.7 in (\ref{e:m2}), we obtain
\begin{equation}\label{e:m3}
\mu\lambda_n\omega_{ii}^2=0,\end{equation} whereby, we find
\begin{equation}\label{e:m4}
\omega_{ii}^2=0.\end{equation}

Now, using  Lemma 3.6, we get
\begin{equation}\label{e:m5}
e_2e_n(\omega_{ii}^n)-e_ne_2(\omega_{ii}^n)=-\omega_{n2}^{1}e_1(\omega_{ii}^n).\end{equation}

Using (\ref{e:q15}), (\ref{e:k50}),  (\ref{e:k51}), (\ref{e:k52}),
(\ref{e:k22}), (\ref{e:b3}) and Lemma 3.7 in (\ref{e:m5}), we obtain
\begin{equation}\label{e:m6}
\omega_{ii}^1=0.\end{equation}

Using (\ref{e:m4}) and (\ref{e:m6}) in (\ref{e:k53}), we find
$\lambda_i=0$ for all $i\in C$. Using this in (\ref{e:q21}), we get
$H=0$, a contradiction. Whereby proof of  Theorem 1.1 is complete.


\bibliography{xbib}

\begin{thebibliography}{99}
\normalsize
\bibitem{r1} A. Arvanitoyeorgos, F. Defever, G.
Kaimakamis, V. Papantoniou. Biharmonic Lorentzian hypersurfaces in
$E_{1}^{4}$, Pac. J. Math. 229(2) (2007), 293-305.

\bibitem{b6} A. Z. Petrov, Einstein spaces, Pergamon
Press, Oxford; 1969.

\bibitem{13} A. Upadhyay, N. C. Turgay, A
classification of biconservative hypersurfaces in a pseudo-Euclidean
space, J. Math. Anal. Appl., 444 (2016), 1703-1720.

\bibitem{biha7} B. Y. Chen, Total mean curvature and submanifolds of finite type, 2nd edition, World Scientific, Hackensack, NJ,
2015.
\bibitem{r3} B. Y. Chen. Classification of
marginally trapped Lorentzian flat surfaces in  $E_{1}^{4}$   and
its application to biharmonic surfaces,  J. Math. Anal. Appl. Vol
340 (2008), 861-875.
\bibitem{r5} B. Y. Chen, S. Ishikawa:
Biharmonic pseudo-Riemannian submanifolds in pseudo-Euclidean
spaces, Kyushu J. Math. 52 (1998), 1-18.
\bibitem{biha6} B. Y. Chen, S.
Ishikawa. Biharmonic surfaces in pseudo-Euclidean spaces. Mem. Fac.
Sci. Kyushu Univ. A. 45 (1991), 323-347.
\bibitem{k1} Deepika, On biconservative Lorentz hypersurface with non-diagonalizable shape operator, Mediterr. J. Math. (2017) 14: 127.
 doi:10.1007/s00009-017-0926-6.
\bibitem{r23} Deepika, Ram Shankar Gupta, A. Sharfuddin. Biharmonic
hypersurfaces with constant scalar curvature in $E^5_s$ , Kyungpook
Math. J., 56(2016), 273-293.
\bibitem{r20} Deepika, Ram Shankar Gupta. Biharmonic
hypersurfaces in $\mathbb{E}^5$ with zero scalar curvature, Afr.
Diaspora J. Math. 18(1) (2015), 12-26.

\bibitem{3} D. Hilbert, Die Grundlagen der Physik, Math. Ann. 92 (1924), 1–32.
\bibitem{r7} F. Defever, G. Kaimakamis, V. Papantoniou.
Biharmonic hypersurfaces of the 4-dimensional semi-Euclidean space
$E_{s}^{4}$, J. Math. Anal. Appl. 315 (2006), 276-286.

\bibitem{4}  G. Y. Jiang, The conservation law for 2-harmonic maps between
Riemannian manifolds, Acta Math. Sin., 30 (1987), 220-225.
\bibitem{r13} I. Dimitri´c.
Submanifolds of  $E^{n}$ with harmonic mean curvature vector, Bull.
Inst. Math. Acad. Sin  20 (1992), 53-65.

\bibitem{5} J. Eells, J. Sampson, Harmonic mappings of Riemannian
manifolds, Amer. J.   Math., 86 (1964), 109-160.
\bibitem{1} J. Marsden, F. Tipler, Maximal
hypersurfaces and foliations of constant mean curvature in general
relativity, Bull. Amer. Phys. Soc., 23 (1978), p. 84.

\bibitem{15} M. A. Magid, Lorentzian isoparametric hypersurfaces", Pacific J.
Math., 118 (1985), 165-197.
\bibitem{k2} N.C. Turgay,  H-hypersurfaces with three distinct principal curvatures in the Euclidean
spaces,  Annali di Matematica, 194 (2015), 1795-1807.
\bibitem{biha15} Ram Shankar Gupta, On biharmonic hypersurfaces in Euclidean
space of arbitrary dimension, Glasgow Math. J. 57 (2015),  633-642.
\bibitem{r14} Ram Shankar Gupta. Biharmonic hypersurfaces in $E_{s}^{5}$, An. St. Univ. Al. I. Cuza, Tomul LXII, 2016, f. 2, vol. 2, 585-593.
\bibitem{r18} Ram Shankar
Gupta. Biharmonic hypersurfaces in $\mathbb{E}^6$ with constant
scalar curvature, International J. Geom., 5(2) 2016, 39-50.
\bibitem{r19} Ram Shankar Gupta, A. Sharfuddin. Biharmonic hypersurfaces in
Euclidean space $\mathbb{E}^5$, J. Geom. 107 (2016), 685-705.
\bibitem{2} S. Stumbles, Hypersurfaces of constant mean
extrinsic curvature, Ann. of Physics, 133 (1981), 28-56.
\bibitem{12} T. Sasahara, Tangentially biharmonic Lagrangian H-umbilical
submanifolds in complex space forms, Abh. Math. Semin. Univ. Hambg.,
85 (2015), 107-123.
\bibitem{6} T.
Hasanis, T. Vlachos, Hypersurfaces in $\mathbb{E}^{4}$ with harmonic
mean curvature  vector field, Math. Nachr. 172 (1995), 145-169.\\
 \end{thebibliography}


Author's address:\\
\textbf{Ram Shankar Gupta}\\
Assistant Professor, University School of Basic and Applied
Sciences, Guru Gobind Singh Indraprastha University, Sector-16C,
Dwarka, New Delhi-110078,
India.\\
\textbf{Email:} ramshankar.gupta@gmail.com\\
\\
\textbf{Prof. Sharfuddin Ahmad}\\Department of Mathematics, Faculty
of Natural Sciences, Jamia Millia Islamia (Central University),
New Delhi-110025, India.\\
\textbf{Email:} sharfuddin\_ahmad12@yahoo.com

\end{document}